\documentclass[12pt,a4paper]{article}

\usepackage{cmap}
\usepackage[margin=20mm,nohead]{geometry}
\usepackage[T1]{fontenc}
\usepackage[utf8]{inputenc}
\usepackage{ae}
\usepackage{graphicx}
\usepackage[american]{babel}
\usepackage{amsfonts}
\usepackage{amsthm}
\usepackage{amsmath}
\usepackage{amssymb}
\usepackage{enumitem}
\usepackage{wasysym}
\usepackage{multirow}

\newtheorem{Def}{Definition}
\newtheorem{Theorem}{Theorem}
\newtheorem{Lemma}[Theorem]{Lemma}
\newtheorem{Corollary}[Theorem]{Corollary}
\newtheorem{Proposition}[Theorem]{Proposition}
\newtheorem{Remark}[Theorem]{Remark}
\newtheorem{Question}{Question}

\newcommand{\A}{\mathcal{A}}
\newcommand{\E}{\mathbb{E}}
\newcommand{\Ex}{\mathcal{E}}
\newcommand{\F}{\mathcal{F}}
\newcommand{\N}{\mathbb{N}}
\newcommand{\R}{\mathbb{R}}
\newcommand{\T}{\mathcal{T}}
\newcommand{\eps}{\varepsilon}
\newcommand{\out}{\operatorname{out}}
\newcommand{\Var}{\operatorname{Var}}
\newcommand{\Cov}{\operatorname{Cov}}
\newcommand{\supp}{\operatorname{supp}}
\newcommand{\bisect}{\operatorname{bisection}}
\newcommand{\BS}{\operatorname{BS}}
\DeclareMathOperator{\LA}{LA}
\DeclareMathOperator{\RLA}{RLA}
\DeclareMathOperator{\MLA}{MLA}
\DeclareMathOperator{\PLA}{PLA}
\DeclareMathOperator{\WLA}{WLA}
\DeclareMathOperator{\PMLA}{PMLA}
\DeclareMathOperator{\MRLA}{MRLA}
\DeclareMathOperator{\PRLA}{PRLA}
\DeclareMathOperator{\WRLA}{WRLA}
\DeclareMathOperator{\PMRLA}{PMRLA}
\DeclareMathOperator{\FF}{fF}
\DeclareMathOperator{\IF}{iF}
\DeclareMathOperator{\FC}{fC}
\DeclareMathOperator{\IC}{iC}

\title{Variants of local algorithms on sparse graphs}
\author{Endre Cs\'oka\\ 
\small{Alfr\'ed R\'enyi Institute of Mathematics}\\
\small{csokaendre@gmail.com}
}
\date{}

\begin{document}

\maketitle

\begin{abstract}
Suppose we want to construct some structure on a bounded-degree graph, e.g., an almost maximum matching, and we want to decide about each edge depending only on its constant-radius neighborhood. 
We examine and compare the strengths of different extensions of these local algorithms. 
A common extension is to use preprocessing, which means that we can make some calculation about the whole graph, and each local decision can also depend on this calculation. 
In this paper, we show that preprocessing is needless: if a nearly optimal local algorithm uses preprocessing, then the same can be achieved by a local algorithm without preprocessing, but with a global randomization.
\end{abstract}



\section{Introduction}

We define a deterministic local algorithm (LA) on graphs or decorated graphs so that we output a decision at each vertex or edge depending on the isomorphism class of its constant-radius neighborhood.
In a random local algorithm (RLA), we assign independent uniform random seeds from $[0, 1]$ to each vertex, and we make a decision at each local structure (typically vertex or edge) depending on its constant-radius seeded neighborhood. 
For example, choosing the set of vertices with higher seeds than all of its neighbors is an RLA that finds an independent set of expected size $\sum\limits_{x \in V} \frac{1}{\deg(x) + 1}$. 

Local algorithms were defined by Linial \cite{Linial} as distributed algorithms using a limited number of synchronized rounds, but otherwise, with no limitation on the computational time and space. 
A distributed algorithm on graphs uses a processor at each vertex, and two processors can directly communicate if they are at neighboring nodes. 
Eventually, each processor makes a decision, which is the output of the algorithm. 
For example, if we want to find a large independent set, then each processor decides whether to put the node into the set or not. 
Several inequivalent variants of distributed algorithms have been studied.

Throughout this paper, local algorithm means a constant-time local algorithm in Linial's sense. 
We study variants of local algorithms on bounded-degree graphs. 
An equivalent definition of local algorithms is that the output of each node is a function of (the isomorphism type of) the constant-radius neighborhood of the node. 

Research on local algorithms was pioneered by Angluin \cite{Angluin}, Linial \cite{Linial}, and Naor and Stockmeyer \cite{NaSt}. 
Angluin \cite{Angluin} studied the limitations of anonymous networks without any unique identifiers. 
Linial \cite{Linial} proved some negative results for the variant where each node has a unique identifier. 
Naor and Stockmeyer \cite{NaSt} presented the first nontrivial positive results. 
For more details about local algorithms, see the survey paper by Suomela \cite{Suomela}. 

For typical problems, we do not expect strictly optimal solutions from local algorithms, but rather approximate solutions. 
For example, we say that we can find an almost maximum independent set if for each $\eps > 0$, there exists a (randomized) local algorithm that outputs an independent set, and the (expected) size of this set is at most $\eps n$ smaller than the size of the maximum independent set. 
In other words, the error is measured in terms of the independence ratio, which is the ratio of the size of the maximum independent set and the size of the vertex set of the graph. 

Local algorithm with preprocessing means that the output of each node is a function of (the isomorphism type of) the graph and the constant-radius neighborhood of the node. 
Equivalently, each vertex receives the same ``global information'' depending on the entire graph, and then they make a constant number of synchronous communication rounds before presenting the output. 
This can also be interpreted as a service:  
there is a center that can provide arbitrary information about the entire graph. 
Taking the maximum independent set problem as an example, each vertex can ask whether it is in the set, and the center should answer using its preprocessed information from the graph and (the isomorphism type of) the constant-radius rooted neighborhood of the node. 
We require that these answers be \emph{consistent}, namely no two neighboring nodes should receive ``yes'', and \emph{nearly optimal}, that is, the proportion of nodes receiving ``yes'' should be close to the relative size of the maximum independent set, that is, the independence ratio. 

In this paper, we show that preprocessing is needless. 
More precisely, if there exists a local algorithm using preprocessing, then there exists another local algorithm with the same radius that the only ``preprocessing'' is a random variable with uniform distribution on $[0,1]$, and it provides an output with at most the same relative error from the optimum in expectation. 
This random variable can also be interpreted as drawing a random sample from a probability distribution on local algorithms, and then applying this local algorithm at each node. 

While this random mixture of local algorithms was not considered prior to this paper, randomness is a powerful and classical technique in the design of distributed algorithms, and particularly useful in breaking the symmetry \cite{ABI, IsIt, Luby}. 
For example, on transitive graphs (such as cycles), any local algorithm should choose the same output at each node. 
Therefore, it is impossible to choose a positive fraction of independent vertices by using an LA. 
However, this is possible with an RLA. 
Our result on preprocessing being irrelevant applies to RLAs, too.

Local algorithms are useful for parameter testing, as well. 
We give a brief overview about parameter testing and the theory of very large graphs; for more details, cf. the survey paper and book \cite{Lovasz, Lovaszbook} by Lov\'asz. 
Parameter testing is an important concept in the theory of bounded-degree graphs~\cite{BeSchSh, BObT, CzuShaSo, Elek2, Elek, LOW, MaRo, NgOn}. 
For a graph parameter, a tester is an algorithm which gets the constant-radius neighborhoods of a constant number of random nodes as input, and outputs a number as an estimation for the parameter. 
We call a parameter testable if for each $\eps > 0$, there exists a tester which estimates the correct value of the parameter with at most $\eps$ error with probability at least $1 - \eps$.

Many of the parameters investigated in the literature come from a maximization problem. 
Examples include the size of the maximum matching, the size of the maximum independent set, or the size of the maximum cut, normalized by the number of nodes. 
Nguyen and Onak \cite{NgOn} proved the testability of several parameters using the following observation. 
If we have a random local algorithm which provides an almost optimal structure (e.g., an almost maximum independent set or matching) then the relative size of the maximum optimal structure is a testable parameter. 
The tester takes the constant-radius neighborhoods of the constant number of random nodes or edges, with the same radius that the random local algorithm uses. 
For each neighborhood, we assign random numbers to the vertices, then we calculate whether the algorithm would put the root vertex or root edge into the structure (e.g., the root vertex into the independent set, or the edge into the matching). 
Then the ratio of these nodes gives a good approximation for the relative size of the optimal structure.

Elek~\cite{Elek} used a random local algorithm with the preprocessing of a finite statistics of constant-radius neighborhoods of random nodes, namely the output at each vertex could depend on its constant-radius neighborhood and this statistics. 
The point of this concept is that the existence of such an algorithm giving good approximation still implies  testability. 
The tester can use the same number of the same radius neighborhoods to make the statistics, and we give this statistics to one further neighborhood, and we calculate the decision at each vertex. 
We repeat this procedure a large number of times, which yields an estimation to the parameter. 
This observation was used as a technique to convert some results in Borel graph theory to theorems in the field of constant-time algorithms, see~\cite{ElekLipi}. 

Summarizing, local algorithms provide a natural tool to solve central problems, and sometimes it has been used with a specific kind of randomization and with some kind of preprocessing. 
We show that for a general class of problems, whether we use randomization or not, preprocessing is needless, because it can be replaced by one public random variable. 

\section{Model and results} \label{model}

Throughout the paper, \textbf{graph} means finite graph with degrees bounded by an absolute constant, allowing loops and parallel edges. 
Whenever we have a function depending on a graph, we mean that it depends on the isomorphism type of the graph. 
In other words, graphs are considered to be unlabelled, unless stated otherwise. 
The $r$-neighborhood of a vertex $x$ of a graph $G$, denoted by $B_r(x)$ or $B_r(G, x)$, means the rooted subgraph of $G$ spanned by all nodes at distance at most $r$ from $x$, and rooted at $x$. 
(Hence, in this case, we label the root, but no other vertex is labelled.)
For a family $\F$ of graphs, denote the family of rooted r-neighborhoods by $\F_r = \big\{B_r(G, x) \big| G \in \F;\ x \in V(G)\big\}$. For any sequence of graphs $G_i$, let $\bigcup G_i$ denote their disjoint union, that is, $V(\bigcup G_i) = \big\{(x, i) \big| x \in G_i\big\}$ and $E(\bigcup G_i) = \Big\{\big((x, i), (y, i)\big) \Big| (x, y) \in \E(G_i)\Big\}$. A 5-tuple $(\F, C, \delta, \A, v)$ is called a \textbf{local choice problem}, where
\begin{itemize}
\item $\F$ is a union-closed family of graphs, namely, $G, H \in \F \Rightarrow G \cup H \in \F$;
\item $C$ is an arbitrary set (the set of choices);
\item $\delta$ is a positive integer (the radius);
\item $\A$ is a set of pairs $(H, c)$ where $H \in \F_{\delta}$ and $c$ is a function $V(H) \rightarrow C$ (the allowed pairs);
\item $v$ is a function $C \rightarrow (-\infty, M]$ for some fixed $M\in \mathbb{R}$ (the value).
\end{itemize}

By \textbf{choice} we mean a function $c: V(G) \rightarrow C$. 
Given a graph $G$, we call a choice $c$ \textbf{allowed} if $\forall x \in V(G): \big(B_{\delta}(G, x), c\upharpoonright_{V(B_{\delta}(G, x))}\big) \in \A$. 
We denote the set of all allowed choices $c$ by $\A(G)$. 
The \textbf{value of a choice} is
\begin{equation} \label{vgdef}
\bar{v}(G, c) = \frac{1}{\big|V(G)\big|} \sum_{x \in V(G)} v\big(c(x)\big),
\end{equation}
and the value of a graph is
\begin{equation} \label{vstardef}
v^*(G) = \sup_{c \in \A(G)} \bar{v}(G, c).
\end{equation}
Given a local choice problem, our goal is to find for any input graph $G$ an allowed choice $c$ with small error $v^*(G)-\bar{v}(G, c)$.

As an illustration, we describe the maximum matching problem in this language. 
Let $\F$ be the family of all graphs, $C = [0, 1] \cup \{\emptyset\}$, and $\delta = 1$. 
Then $(H, c) \in \A$ iff for the root $x$ of $H$, $c(x) = \emptyset$ or there exists exactly 1 neighbor $y$ of $x$ with $c(x) = c(y)$. 
Finally, $v(col)$ is 0 if $col = \emptyset$ and $v(col) = \frac{1}{2}$ otherwise. 
Then the allowed choices $c$ describe the matchings: $c(x) = \emptyset$ if $x$ is unmatched, otherwise $x$ is matched with the unique neighboring vertex $y$ such that $c(x) = c(y)$. 
The value $\bar{v}(G, c)$ is the relative size of this matching, normalized by $\big|V(G)\big|$.

At first, it might seem to be more natural and more general to define $v: \A \rightarrow (-\infty, M]$ and $\bar{v}(G, c) = \frac{1}{|V(G)|} \sum\limits_{x \in V(G)} v\big(B_{\delta}(G, x), c\upharpoonright_{B_{\delta}(G, x)}\big)$. 
In fact, this definition would not be more general than the original version. 
Roughly speaking, we can define the coloring so as to include the value of the coloring at the point. 
More formally, let $(\F, C, \delta, \A, v)$ be an extended local choice problem, namely we have this more general $v$. 
Let $C' = C \times \R$ and $\A' = \Big\{ \big(H, (c_1, c_2)\big) \Big| (H, c_1) \in \A;\ c_2\big(root(H)\big) = v(H, c_1) \Big\}$, and $v'\big(H,(c_1, c_2)\big) = c_2$. 
Then the local choice problem $(\F, C', \delta, \A', v')$ is equivalent to the extended local choice problem $(\F, C, \delta, \A, v)$. 
The details are left to the Reader.

Now we define different versions of local algorithms for finding such an allowed choice $c$. 
We assign independent identically distributed random variables to the vertices with a fixed distribution $D$. 
We denote this random assignment by $\omega: V(G) \rightarrow \Omega$. 
The most important case is when $D$ is a continuous distribution, say, uniform on $[0, 1]$, but it can be a constant number of random bits, or any other distribution. 
We take one more independent public random variable $g$ with uniform distribution on $[0, 1]$. 
The choice $g \sim U[0,1]$ is without loss of generality, because every probability distribution can be written as a function of $g$. 
Alternatively, we could simply allow that $g$ be drawn from any distribution. 

\begin{Def}{}
[\textbf{Preprocessed}] [\textbf{Mixed}] [\textbf{Random}] \textbf{Local Algorithm} ([\textbf{P}][\textbf{M}][\textbf{R}]\textbf{LA}). For a fixed radius $r$, we set $c(x)$ depending on $B_r(x)$
\begin{itemize}
\item and $G$ if the local algorithm is Preprocessed. 
\item and $g$ if the local algorithm is Mixed. 
\item and $\omega|_{V(B_r(x))}$ if the local algorithm is Random.
\end{itemize}
The choice made by the $[P][M][R]\LA$ $l$, given the graph $G$, the random vector $\omega$, and the public random variable $g$ is denoted by $l[G, \omega, g]$. 
\end{Def}

Now we have $2\times 2\times 2$ different types of local algorithms. 
We say that an algorithm is correct if it always produces allowed choices. 
(We could use ``with probability 1'' instead of ``always'' with essentially the same proofs.)
The main result of this paper is that $\MRLA$ and $\PMRLA$ are equally strong, in the following sense. 
We say that a local choice problem is approximable in one of the eight types of local algorithms defined above, if for all $\eps > 0$, there exists a correct local algorithm $l$ in that type such that $\forall G \in \F: \E\big(\bar{v}(G, l[G, \omega, g])\big) \ge v^*(G) - \eps$.

\begin{Theorem} \label{main}
Let $\F$ be a union-closed family of finite graphs with a given degree bound. 
If a local choice problem for graphs in $\F$ is approximable in $\PMRLA$, then it is also approximable in $\MRLA$.
\end{Theorem}

\begin{Remark}
As an immediate corollary, we obtain that the optimal value corresponding to a local choice problem that is approximable in $\PMRLA$ is testable. 
Indeed, we can switch to an $l'\in \MRLA$ according to Theorem~\ref{main}, and the expectation of the output $\E\big(\bar{v}(G, l[\omega, g])\big)$ is a testable parameter of the graph $G$. 
\end{Remark}

\begin{figure}
\begin{center}
\includegraphics[scale=0.6]{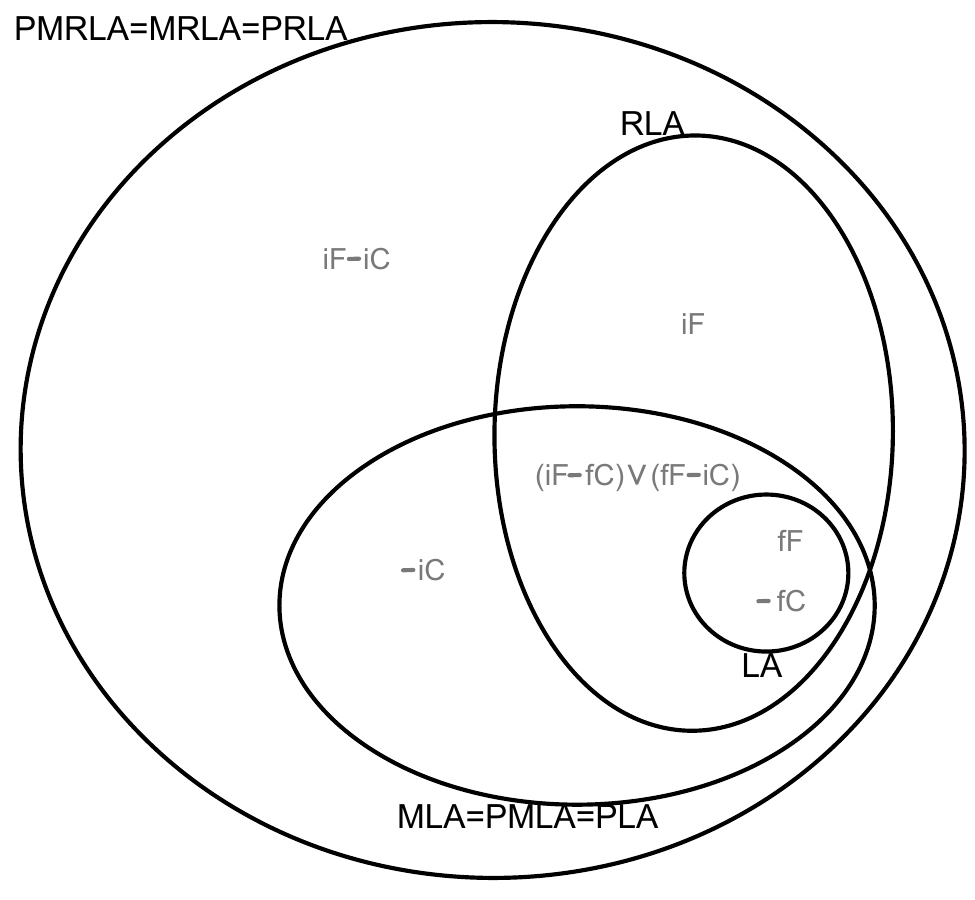}
\end{center}
\caption{Hierarchy of local choice problems, where $\F$ consists of all graphs with a given degree bound.}
\label{fig:Venn}
\end{figure}

In Section~\ref{sec:gen}, we present a more general form of the main result Theorem~\ref{main}. 
In Section~\ref{sec:collapse}, we are going to show that the only collapses in the hierarchy of local algorithms are those implied by $M=P=PM$; cf. Figure~\ref{fig:Venn}.

\section{Proof of the main results}\label{sec:gen}

\begin{Theorem} \label{general}
Let $\F$ be a union-closed family of finite graphs with a given degree bound. 
Let $\eps>0$ and let $b: \mathbb{R} \rightarrow \mathbb{R}$ be a strictly monotone increasing concave function. 
Assume that there exists a correct $l\in \PMRLA$ such that for each graph $G \in \F$ we have $\E_{\omega, g}\Big(\bar{v}\big(G, l[G, \omega, g]\big)\Big) \ge b\big(v^*(G)\big)$. 
Then there exists a correct $l'\in \MRLA$ using the same radius $r$ and the same distribution of $\omega$ such that $\E_{\omega, g}\Big(\bar{v}\big(G, l'[\omega, g]\big)\Big) > b\big(v^*(G)\big) - \eps$.
\end{Theorem}

Recall that the distribution $D$ of $\omega$ is arbitrary but fixed. 
Therefore, [P][M] random local algorithms can have different strengths depending on $D$, but we always use the same $D$. 
For example, by applying the two theorems when $D$ is a constant distribution, we obtain the corollary that $\PMLA$ and $\MLA$ have the same power.

At the end of the paper, we will show that preprocessing or mixing or both are just equally strong, but useful and not exchangeable with randomizing in any direction. For example, we show a problem which is approximable by $\MLA$ but not by $\RLA$. 

Finally, we note that our results extend naturally to bounded-degree uniform random graphs (URGs). 
An URG is an involution invariant probability distribution $\sigma$ on countable rooted graphs with a given degree bound. 
Involution invariance means that by drawing a sample from $\sigma^*$, which is the normalized version of $\sigma$ by the degree of the root, and then moving the root to a random neighbor, we obtain the same distribution $\sigma^*$ on the set of bounded-degree, countable rooted graphs. 
Every finite graph can be viewed as an URG by choosing $\sigma$ as the uniform distribution on the vertices as roots. 
For more details, cf. \cite{Lovaszbook}. 
We mention the necessary modification in the assertions and proofs. 
By making these modifications, all arguments can be applied directly to URGs. 
Local algorithms and valuations generalize trivially to URGs: the output depends on the constant-radius neighborhood of the (random) root, and the value of a choice is $\bar{v}(G, c) = \int v\big(c(x)\big)\, d\sigma(x)$, where $x$ is the root of the rooted graph sampled from $\sigma$. 
Note that by considering a finite graph an URG as above, this formula translates to the original definition of $\bar{v}(G, c)$. 
Instead of a union-closed family of graphs $\F$, we need that the family $\F$ of URGs be convex. 
In fact, it would be enough to assume that $\F$ is closed under (finite) convex combinations with rational coefficients. 
Everything else carries from finite graphs to URGs verbatim. 

\begin{Theorem} \label{generaluni}
Let $\F$ be a convex set of unimodular random graphs with a given degree bound. 
Let $\eps>0$ and let $b: \mathbb{R} \rightarrow \mathbb{R}$ be a strictly monotone increasing concave function. 
Assume that there exists a correct $l\in \PMRLA$ such that for each graph $\sigma \in \F$ we have $\E_{\omega, g}\Big(\bar{v}\big(\sigma, l[\sigma, \omega, g]\big)\Big) \ge b\big(v^*(\sigma)\big)$. 
Then there exists a correct $l'\in \MRLA$ using the same radius $r$ and the same distribution of $\omega$ such that $\E_{\omega, g}\Big(\bar{v}\big(\sigma, l'[\omega, g]\big)\Big) > b\big(v^*(\sigma)\big) - \eps$. 
In particular, if a local choice problem for URGs in $\F$ is approximable in $\PMRLA$, then it is also approximable in $\MRLA$.
\end{Theorem}

\subsection{Proof of Theorems~\ref{main} and \ref{general}}

We assume that $r \ge \delta$. 
Denote by $s_r(G)$ the distribution of $B_r(G, x)$ for a uniform random vertex $x \in V(G)$, namely $s_r(G)(H) = \Big|\big\{x \in V(G) \mid B_r(G, x) \cong H\big\}\Big| / \big|V(G)\big|$ for any $H \in \F_r$. 
For URGs, $s_r(G)$ denotes the distribution of $B_r(G, x)$ for a random rooted graph $G$ drawn from $\sigma$ with root $x$. 

\begin{Lemma} \label{linear}
Given an $l\in \MRLA$ using radius $r$, the expectation $\E_{\omega, g}\Big(\bar{v}\big(G, l[\omega, g]\big)\Big)$ is a linear function of $s_r(G)$. 
\end{Lemma}

\begin{proof}
The function $v$ is bounded from above by $M$, thus the expected value exists. 
We note that it might be $-\infty$; however, when we apply the lemma in the proof of the main theorems, this can never be the case. 
Given an $x\in V(G)$, the value $v\big(l[\omega, g](x)\big)$ only depends on $B_r(G, x)$, $\omega$, and $g$. 
Consequently, $\E_{\omega, g}\Big(v\big(l[\omega, g](x)\big)\Big)$ only depends on $B_r(G, x)$. 
Define 
\begin{equation} \label{pdef}
p_l\big(B_r(G, x)\big) = \E_{\omega, g}\Big(v\big(l[\omega, g](x)\big)\Big).
\end{equation}
Then 
\begin{equation*}
\E_{\omega, g}\Big(\bar{v}\big(G, l[\omega, g]\big)\Big) \mathop{=}^{\eqref{vgdef}} \E_{\omega, g}\Big(\frac{1}{\big|V(G)\big|} \sum_{x \in V(G)} v\big(l[\omega, g](x)\big)\Big)
\end{equation*}
\begin{equation*} 
= \frac{1}{\big|V(G)\big|} \sum_{x \in V(G)} \E_{\omega, g}\Big(v\big(l[\omega, g](x)\big)\Big) = \frac{1}{\big|V(G)\big|} \sum_{x \in V(G)} p_l\big(B_r(G, x)\big) = \sum_{H \in \F_r} s_r(G)(H) \cdot p_l(H).
\end{equation*}
\end{proof}

We first prove a quite technical lemma.

\begin{Lemma} \label{correct}
For each local choice problem and radius $r$, there exists a graph $T_r \in \F$ such that for all $G \in \F$, whenever an $l\in \MRLA$ with radius $r$ produces an allowed choice on $T_r \cup G$, then $l$ is correct. 
Moreover, $s_r(T_r)(H)>0$ for every $H\in \F_r$. 
\end{Lemma}

\begin{proof}
For each $H \in \F_{\delta+r}$, let us choose a graph $a(H)$ so that $\exists x \in V\big(a(H)\big): B_{\delta +r}\big(a(H), x\big) \cong H$. We show that $T_r = \bigcup\limits_{H \in \F_{\delta+r}} a(H)$ satisfies the requirement.

Suppose that an $l\in \MRLA$ is not correct. 
This means that there exist $G \in \F$, $x \in V(G)$, $\omega: V(G) \rightarrow \supp(\Omega)$, and $g$ such that $\Big(B_{\delta}(x), l[\omega, g]\upharpoonright_{V(B_{\delta}(x))}\Big) \notin \A$. 
For each $y \in V\big(B_{\delta}(x)\big)$, $l[\omega, g](y)$ only depends on $B_{\delta + r}(y)$, $\omega\upharpoonright_{V(B_{\delta + r}(y))}$, and $g$. 
Since $V\big(B_r(y)\big) \subseteq V\big(B_{\delta + r}(x)\big)$, the tuple $\Big(B_{\delta}(x), l[\omega, g]\upharpoonright_{V(B_{\delta}(x))}\Big)$ only depends on $B_{\delta + r}(x)$, $\omega\upharpoonright_{V(B_{\delta + r}(x))}$ and $g$. 
Hence, if we take the component $a\big(B_{\delta + r}(x)\big)$ of $T_r$, the same $\omega$ on $B_{\delta + r}(x')$ where $x'$ is the vertex in $T_r$ corresponding to $x$ in $B_{\delta + r}(x)$), and the same $g$, then it produces the same pair

\noindent
$\Big(B_{\delta}(T_r, x'), l[\omega, g]\upharpoonright_{V(B_{\delta}(T_r, x'))}\Big) \cong \Big(B_{\delta}(x), l[\omega, g]\upharpoonright_{V(B_{\delta}(G, x))}\Big) \notin \A$.
\end{proof}

Before proving the next lemma, observe that
\begin{equation*}
v^*(G \cup H) = \frac{|V(G)| v^*(G) + |V(H)| v^*(H)}{|V(G)| + |V(H)|}\text{ and }s_r(G \cup H) = \frac{|V(G)| s_r(G) + |V(H)| s_r(H)}{|V(G)| + |V(H)|}.
\end{equation*}
Also note that we use the notion of $\limsup$ for a function somewhat unconventionally. 
Usually, even if $x$ is in the domain of the function $f$, the expression $\limsup\limits_{x_n\rightarrow x} f(x_n)$ means that we take the limit superior over all sequences in the pointed neighborhood of $x$, that is, $x_n\neq x$ is required. 
However, we allow $x_n=x$ if $x$ is in the domain, and in particular, $\limsup\limits_{x_n\rightarrow x} f(x_n) \geq f(x)$. 
In fact, because we apply the notion to a concave function whose domain is closed under finite convex combinations with rational coefficients, this distinction should not cause any difference in the result. 
However, we prefer to avoid such unnecessary technical complications. 

\begin{Lemma}\label{lem:concon}
Given $\F$ and $r$, let $T_r$ be as in Lemma~\ref{correct}, and let $\tilde\F=\{G\cup T_r\mid G\in \F\}$. 
Let $S_r = cl \big( \big\{ s_r(G) \mid G \in \F \big\} \big)$, $\tilde{S}_r = cl \big( \big\{ s_r(G) \mid G \in \tilde{\F} \big\} \big)$, and define the function $m_r: \tilde{S}_r \rightarrow \R$ by 
\begin{equation} \label{mdef}
m_r(q) = \limsup\limits_{\substack{s_r(G_n) \rightarrow q\\ G_n\in \tilde{\F}}} v^*(G_n).
\end{equation}
Then $S_r$ is convex, $S_r = \tilde{S}_r$, and $m_r$ is concave and continuous.
\end{Lemma}

\begin{proof}
For a $k\in \mathbb{N}$ and a graph $G$, let $k \times G$ denote $\bigcup\limits_{i = 1}\limits^k G_i$, where each $G_i$ is isomorphic to $G$. 
As $\lim\limits_{k \rightarrow \infty} s_r((k \times G) \cup T_r) = s_r(G)$, we have $S_r = \tilde{S}_r$. 

Given the choices $c_i: V(G_i) \rightarrow C$, let $\sum\limits_{i=1}\limits^k c_i: \bigcup\limits_{i=1}\limits^k V(G_i) \rightarrow C$ denote the function defined by $(\sum\limits_{i=1}\limits^k c_i)\big((x, j)\big) = c_j(x)$. For a choice $c: V(G) \rightarrow C$, let $k \times c = \sum\limits_{i=1}\limits^k c_i$, where each $c_i: G_i \rightarrow C$ is a copy of $c: G \rightarrow C$.

Let $q_0, q_1 \in S_r$, and for all $\lambda \in [0, 1]$ define $q_{\lambda} = (1 - \lambda) \cdot q_0 + \lambda \cdot q_1$. 
Then

\begin{equation*}
(1 - \lambda) \cdot m(q_0) + \lambda \cdot m(q_1) \mathop{=}^{\eqref{mdef}} (1 - \lambda) \cdot \limsup\limits_{\substack{s_r(G_n) \rightarrow q_0\\ G_n\in \tilde{\F}}} v^*(G_n) + \lambda \cdot \limsup\limits_{\substack{s_r(G_n) \rightarrow q_1\\ G_n\in \tilde{\F}}} v^*(G_n)
\end{equation*}
\begin{equation*}
\mathop{=}^{\eqref{vstardef}} (1 - \lambda) \cdot \limsup\limits_{\substack{s_r(G_n) \rightarrow q_0\\ G_n\in \tilde{\F}}} \sup_{c \in \A(G_n)} \bar{v}(G_n, c) + \lambda \cdot \limsup\limits_{\substack{s_r(G_n) \rightarrow q_1\\ G_n\in \tilde{\F}}} \sup_{c \in \A(G_n)} \bar{v}(G_n, c)
\end{equation*}
\begin{equation*}
= \limsup \Big\{ (1 - \lambda) \cdot \bar{v}(G^{(0)}_n, c^{(0)}_n) + \lambda \cdot \bar{v}(G^{(1)}_n, c^{(1)}_n) \Big| \forall i \in \{0, 1\}: \big(s_r(G^{(i)}_n) \rightarrow q_i;\ c^{(i)}_n \in \A(G^{(i)}_n) \big)\Big\}
\end{equation*}
\begin{equation*}
= \limsup \Big\{ \frac{b_n - a_n}{b_n} \cdot \bar{v}(G^{(0)}_n, c^{(0)}_n) + \frac{a_n}{b_n} \cdot \bar{v}(G^{(1)}_n, c^{(1)}_n)
\end{equation*}
\begin{equation*}
\Big| a_n, b_n \in \N;\ \frac{a_n}{b_n} \rightarrow \lambda;\ \forall i \in \{0, 1\}: \big( s_r(G^{(i)}_n) \rightarrow q_i;\ c^{(i)}_n \in \A(G^{(i)}_n) \big) \Big\}
\end{equation*}
\begin{equation*}
= \limsup \bigg\{ \bar{v}\Big((b_n - a_n)\big|V(G^{(1)})\big| \times G^{(0)}_n \bigcup a_n\big|V(G^{(0)})\big| \times G^{(1)}_n, (b_n - a_n)\big|V(G^{(1)})\big| \times c^{(0)}_n
\end{equation*}
\begin{equation*}
+ a_n\big|V(G^{(0)})\big| \times c^{(1)}_n\Big) \bigg| a_n, b_n \in \N;\ \frac{a_n}{b_n} \rightarrow \lambda;\ \forall i \in \{0, 1\}: \big( s_r(G^{(i)}_n) \rightarrow q_i;\ c^{(i)}_n \in \A(G^{(i)}_n) \big) \bigg\}
\end{equation*}
(We always assume $G_n^{(i)}\in \tilde{\F}$.)
The neighborhood distribution of the graph appearing in the last line is $s_r\Big((b_n - a_n)\big|V(G^{(1)})\big| \times G^{(0)}_n + a_n\big|V(G^{(0)})\big| \times G^{(1)}_n)\Big) =  q_{a_n/b_n} \rightarrow q_{\lambda}$. 
This implies the convexity of $S_r$, and continuing the calculations,
\begin{equation*}
\le \limsup\big\{ \bar{v}(G_n, c_n) \big| s_r(G_n) \rightarrow q_{\lambda};\ c_n \in \A(G_n) \big\} \mathop{=}^{\eqref{vstardef}} \limsup\limits_{\substack{s_r(G_n) \rightarrow q_\lambda\\ G_n\in \tilde{\F}}} v^*(G_n) \mathop{=}^{\eqref{mdef}} m(q_{\lambda}),
\end{equation*}
which verifies the concavity of $m$.

In particular, $m_r$ is lower semicontinuous. (It is not necessarily continuous on the boundary.) 
We show that $m_r$ is upper semicontinuous, as well.

Suppose that $q_n \rightarrow q$. 
According to \eqref{mdef}, for each $n \in \N$, there exists a $G_n \in \tilde{\F}$ so that $\big\| s_r(G_n) - q_n \big\| < \frac{1}{n}$ and $m_r(q_n) - v^*(G_n) < \frac{1}{n}$. 
Then $\lim\limits_{n \rightarrow \infty} s_r(G_n) = q$, and $\limsup\limits_{n \rightarrow \infty} m_r(q_n) = \limsup\limits_{n \rightarrow \infty} v^*(G_n) \le \limsup\limits_{s_r(G_n) \rightarrow q} v^*(G_n) = m_r(q)$, verifying the upper semicontinuity of $m_r$.
\end{proof}

When generalizing the results to URGs, a little care is needed at this point. 
We defined $\tilde\F=\{G\cup T_r\mid G\in \F\}$ for finite graphs, and the point was that in these graphs, every rooted $r$-neighborhood occurs with positive ratio (that occurs with positive ratio in some graph in $\F$). 
For URGs, we should first identify those rooted $r$-neighborhoods that occur as the $r$-neighborhood in an URG $\sigma\in \F$ for a positive measure set of the rooted graphs with respect to $\sigma$. 
Then $\tilde\F$ consists of those URGs $\sigma\in \F$ where all these the $r$-neighborhoods occur with positive probability with respect to $\sigma$. 
As $\F$ is convex, this set contains all inner points of $\F$. 

Before proving the main theorems, we need to make some observations on real valued convex or concave functions defined on convex domains. 

\begin{Lemma} \label{convex2}
Let $Q \subset \R^n$ be a compact convex set, and let $f_0, f_1: Q \rightarrow \R$ be two convex functions such that $\forall q \in Q: f_0(q) > 0 \,\,\text{or}\,\, f_1(q) > 0$. 
Then there is a convex combination of the functions which is positive on $Q$. 
Formally, $\exists \lambda \in [0, 1]: \forall q \in Q: f_{\lambda}(q) = \big((1 - \lambda) \cdot f_0 + \lambda \cdot f_1\big)(q) > 0$.
\end{Lemma}

\begin{proof}
Let $f_{\lambda}^- = \big\{q \in Q \big| f_{\lambda}(q) \le 0\big\}$. 
The sets $f_{\lambda}^-$ are convex and compact, and $f_0^-$ and $f_1^-$ are disjoint. 
If $f_0(q) > 0$ and $f_1(q) > 0$, then $f_{\lambda}(q) > 0$ as well, so $f_{\lambda}^- \subseteq f_0^- \cup f_1^-$. 
Hence, $f_{\lambda}^- \subseteq f_0^-$ or $f_{\lambda}^- \subseteq f_1^-$ for all $\lambda\in[0,1]$.

The function $\lambda \rightarrow \min\limits_{q \in Q} f_{\lambda}(q)$ is continuous, so $\big\{ \lambda \in [0, 1] \big| \min\limits_{q \in Q} f_{\lambda}(q) \le 0 \big\}$ is closed. 
Therefore the sets $A = \big\{\lambda \in [0, 1] \big| f_{\lambda}^- \cap f_0^- \neq \emptyset\big\}$ and $B = \big\{\lambda \in [0, 1] \big| f_{\lambda}^- \cap f_1^- \neq \emptyset\big\}$ are closed, disjoint and nonempty. 
Thus $A \cup B \neq [0, 1]$, because $[0, 1]$ is a connected topological space. 
Therefore, there exists a $\lambda \in [0, 1]\setminus (A\cup B)$, and such a $\lambda$ satisfies the requirements of the lemma.
\end{proof}

\begin{Lemma} \label{convexn}
Let $Q \subset \R^n$ be a compact convex set. For each $q \in Q$ let $f_q: Q \rightarrow \R$ be a convex function such that $f_q(q) > 0$. 
Then there is a finite convex combination of the functions $f_q$ that is positive on $Q$. 
Formally, $\exists k\in\mathbb{N}, q_1, q_2, \ldots, q_k \in Q$ and $\lambda_1, \lambda_2, \ldots, \lambda_k \ge 0$ such that $\sum\limits_i \lambda_i = 1$ and $\forall p \in Q: \sum\limits_i \lambda_i f_{q_i}(p) > 0$. 
\end{Lemma}

\begin{proof}
Consider the set $\T$ of finite convex combinations of the $f_q$. 
Each function in this set is convex. 
For a function $f \in \T$, let $f^+ = \big\{q \in Q \big| h(q) > 0\big\}$ and $f^- = \big\{q \in Q \big| h(q) \le 0\big\}$ denote the positive and the nonpositive set of $f$, respectively. 
The positive sets $f^+$ are open, and they cover the compact set $Q$. 
Thus there exists a finite list of functions $f \in \T$ such that their positive sets cover $Q$. 
Let $f_1, f_2, \ldots, f_h$ be a shortest list.

Assume that $h > 1$, and let $Q' = f_3^- \cap f_4^- \cap \cdots \cap f_h^-$. The set $Q'$ is the intersection of finitely many convex compact sets, so $Q'$ itself is convex and compact. 
At each point $q \in Q'$, $f_1(q) > 0$ or $f_2(q) > 0$, as otherwise $q$ would not be covered by any of the $h_i^+$. 
According to Lemma~\ref{convex2} there exists a convex combination $f_0$ of $f_1$ and $f_2$ which is positive on $Q'$. 
Clearly, $f_0 \in \T$ and $f_0^+ \cup f_3^+ \cup f_4^+ \cup \cdots \cup f_h^+ = f_0^+ \cup (Q - Q') = Q$, contradicting the minimality of $h$. 
Thus $h = 1$, and the one-element list $f$ satisfies the requirements of the lemma.
\end{proof}

\begin{Lemma} \label{linearn}
Let $Q \subset \R^n$ be a compact convex set, and let $f: Q \rightarrow \R$ be concave function. 
For each $q \in Q$, let $f_q: Q \rightarrow \R$ be a linear function such that $f_q(q) > f(q)$. 
Then there exists a finite convex combination of the $f_q$ that dominates $f$. 
Formally, $\exists k\in\mathbb{N}, q_1, q_2, \ldots, q_k \in Q$ and $\lambda_1, \lambda_2, \ldots, \lambda_k \ge 0$ such that $\sum\limits_i \lambda_i = 1$ and $\forall p \in Q:$
\begin{equation} \label{conveq}
\sum_i \lambda_i f_{q_i}(p) > f(p).
\end{equation}
\end{Lemma}

\begin{proof}
The functions $f_q - f$ are convex and positive at $q$, thus Lemma~\ref{convexn} applies.
\end{proof}

We are ready to prove the main results. 

\begin{proof}[\emph{\textbf{Proof of Theorem~\ref{general}}}]
We use the notation introduced in Lemmas~\ref{linear}, \ref{correct}, and \ref{lem:concon}. 
For an algorithm $l\in \MRLA$ and a distribution $q \in S_r$ let 
\begin{equation} \label{udef}
u(l, q) = \sum\limits_{H \in \F_r} q(H) \cdot p_l(H).
\end{equation}
where $p_l(H)$ is as in Lemma~\ref{linear}. 
As $v \le M$, we have $p_l \le M$ and $u \le M$. 

Given a finite graph $\tilde{G}\in \tilde{\F}$, the $\PMRLA$ $l$ simplifies to an $l_{\tilde{G}}\in \MRLA$: we run the algorithm $l$ on any graph $G\in \F$, and whenever a global parameter of $G$ is queried, we provide the parameter value in $\tilde{G}$. 
According to Lemma \ref{correct}, this simplified $\MRLA$ $l_{\tilde{G}}$ is correct for all graphs $G\in \F$. 
By the condition $\E_{\omega, g}\Big(\bar{v}\big(\tilde{G}, l[\tilde{G}, \omega, g]\big)\Big) \ge b\big(v^*(\tilde{G})\big)$ of Theorem~\ref{general}, we have that the value $\E_{\omega, g}\Big(\bar{v}\big(\tilde{G}, l[\tilde{G}, \omega, g]\big)\Big)$ is a proper (i.e., finite) real number. 
Since $s_r(\tilde{G})(H)>0$ for each $H\in \F_r$ in this case, we have that the linear function $u(l_{\tilde{G}}, q)$ is also proper in the sense that all its coefficients $p_{l_{\tilde{G}}}(H)$ are real numbers, rather than $-\infty$, as 
\begin{equation*}
\sum\limits_{H \in \F_r} s_r(\tilde{G})(H) \cdot p_{l_{\tilde{G}}}(H) = u(l_{\tilde{G}}, s_r(\tilde{G}))=\E_{\omega, g}\Big(\bar{v}\big(\tilde{G}, l_{\tilde{G}}[\omega, g]\big)\Big) = \E_{\omega, g}\Big(\bar{v}\big(\tilde{G}, l[\tilde{G}, \omega, g]\big)\Big)\neq -\infty.
\end{equation*} 
The $\MRLA$ $l_{\tilde{G}}$ corresponding to $\tilde{G}$ then yields a choice process as above for any $G\in \F$, and in particular a value $\bar{v}\big(G, l_{\tilde{G}}[\omega, g]\big)$ whose expectation is a linear function according to Lemma~\ref{linear}, with a compact domain. 
Hence, this linear function has an absolute lower bound $K_0\in \mathbb{R}$. 
Thus by putting $K_1=\max(|M|, |K_0|)$, we have $|v^*(G)|\leq K_1$ for any $G\in \F$, and then by the defining equation \eqref{mdef}, we have $|m_r(q)|\leq K_1$ for all $q\in S_r$. 
Let $K=\max(K_1, |b(-K_1)|, |b(K_1)|)$. 
Then $|b(m_r(q))|\leq K$ for all $q\in S_r$.

The function $b$ is continuous, and consequently, it is uniformly continuous on $[-K_1, K_1]$. 
Let $\delta_1>0$ be such that whenever $|x-x'|\leq\delta_1$ for $x,x'\in [-K_1, K_1]$ then $|b(x)-b(x')|\leq\eps/4$. 
Similarly, the function $m_r$ is uniformly continuous on $S_r$, thus there is a $\delta_2>0$ such that whenever $||q-q'||\leq\delta_2$ then $|m_r(q)-m_r(q')|\leq \delta_1$; here $||q-q'||$ is the usual $L_2$ norm of the vector $q-q'$, that is, the Euclidean distance of $q$ and $q'$. 
Finally, let $\delta_0=\min\left(\delta_2/D,\frac{\eps/(8K)}{1+\eps/(8K)}\right)$, where $D$ is the diameter of $S_r$. 

Given a $q\in S_r$, let $U$ be the homothetic image of $S_r$ with center $q$ and ratio $\delta_0$. 
In particular, the diameter of $U$ is smaller than $\delta_2$. 
Assume indirectly that for any $G\in \tilde\F$ with $s_r(G)\in U$ we have $v^*(G)<m_r(s_r(G))-\delta_2$. 
Then in an inner point of $U$ of the form $s_r(G_0)$ with $G_0\in\tilde\F$, the function value of $m_r$ is obtained as the limit of a sequence $v^*(G_n)$, where $G_n\in \tilde\F$, $s_r(G_n)\in U$, and $s_r(G_n)\rightarrow s_r(G_0)$. 
By the continuity of $m_r$, we have $\lim m_r(s_r(G_n)) = m_r(s_r(G_0))$. 
On the other hand, $m_r(s_r(G_0)) = \lim v^*(G_n) \leq \lim m_r(s_r(G_n)) - \delta_2 = m_r(s_r(G_0)) - \delta_2$ by the indirect assumption, a contradiction. 
Thus for the given $q\in S_r$ there exists a $\tilde{G}\in \tilde\F$ such that $s_r(\tilde{G})\in U$ and $v^*(\tilde{G})\geq m_r(s_r(\tilde{G}))-\delta_2$, that is, $|m_r(s_r(\tilde{G})) - v^*(\tilde{G})| \leq \delta_2$. 
Pick such a graph $\tilde{G}$ for the given $q$.  
Then 
\begin{equation*}
b(m_r(q))-u(l_{\tilde{G}}, q) =  
\end{equation*}
\begin{equation*}
\big(b(m_r(q)) - b(m_r(s_r(\tilde{G})))\big) + \big(b(m_r(s_r(\tilde{G}))) - b(v^*(\tilde{G}))\big) + 
\end{equation*}
\begin{equation*}
\big(b(v^*(\tilde{G})) - u(l_{\tilde{G}}, s_r(\tilde{G}))\big) + \big(u(l_{\tilde{G}}, s_r(\tilde{G})) -u(l_{\tilde{G}}, q)\big)\leq
\end{equation*}
\begin{equation}\label{eq:threeterms}
\big|b(m_r(q)) - b(m_r(s_r(\tilde{G})))\big| + \big|b(m_r(s_r(\tilde{G}))) - b(v^*(\tilde{G}))\big| + \big(u(l_{\tilde{G}}, s_r(\tilde{G})) -u(l_{\tilde{G}}, q)\big).
\end{equation}
The estimation is provided by the condition of the theorem: the third summand after the equality sign is $b(v^*(G)) - u(l_{G}, s_r(G)) = b(v^*(G)) - \E_{\omega, g}\Big(\bar{v}\big(G, l[G, \omega, g]\big)\Big)\leq 0$. 
In addition, the first two summands were estimated by their absolute value. 
We show that all three summands in \eqref{eq:threeterms} are at most $\eps/4$. 

As $q$ and $s_r(\tilde{G})$ are both in $U$, a set with diameter at most $\delta_2$, we have $||q-s_r(\tilde{G})||\leq \delta_2$, implying $|m_r(q)-m_r(s_r(\tilde{G}))|\leq \delta_1$, and then $\big|b(m_r(q)) - b(m_r(s_r(\tilde{G})))\big|\leq \eps/4$ by the choice of $\delta_1$ and $\delta_2$. 

We chose $\tilde{G}$ so that $|m_r(s_r(\tilde{G})) - v^*(\tilde{G})| \leq \delta_1$. 
Hence, by the choice of $\delta_1$, the second summand is $\big|b(m_r(s_r(\tilde{G}))) - b(v^*(\tilde{G}))\big|\leq \eps/4$. 

Finally, we estimate the third summand $\big(u(l_{\tilde{G}}, s_r(\tilde{G})) -u(l_{\tilde{G}}, q)\big)$. 
Let $q'$ be the preimage of $s_r(\tilde{G})$ under the above homothecy with center $q$ and ratio $\delta_0$. 
As $\delta_0\leq \frac{\eps/(8K)}{1+\eps/(8K)}$, we have $||q - s_r(\tilde{G})||\leq \eps/(8K)\cdot ||s_r(\tilde{G}) - q'||$. 
Denote the restriction of the linear function $u(l_{\tilde{G}}, .)$ to the segment $[q,q']$ by $u_0$. 
(These are proper linear functions with real coefficients, as $\tilde{G}\in \tilde\F$.) 
Then 
\begin{equation*}
u_0(s_r(\tilde{G})) = u(l_{\tilde{G}}, s_r(\tilde{G})) = \E_{\omega, g}\Big(\bar{v}\big(\tilde{G}, l[\tilde{G}, \omega, g]\big)\Big) \geq b(v^*(\tilde{G}))\geq -K
\end{equation*} 
by the assumption of the theorem and the facts that $b(x)\in [-K, K]$ for all $x\in [-K_1, K_1]$ and $v^*(\tilde{G})\in [-K_1, K_1]$. 
Furthermore, we show that $u_0(p)\leq K$ for all $p\in [q,q']$. 
To this end, it is enough to prove that $u(l_{\tilde{G}}, s_r(G))\leq K$ for all $G\in \F$, since the points $s_r(G)$ form a dense subset in $S_r$ and the function $u(l_{\tilde{G}}, .)$ is linear. 
Indeed, $u(l_{\tilde{G}}, s_r(G)) = \E_{\omega, g}\Big(\bar{v}\big(G, l[\tilde{G}, \omega, g]\big)\Big) \leq M \leq K$. 
In particular, $u_0(q')\leq K$. 
In summary, we have the linear function $u_0$ defined on the segment $[q,q']$, and the point of this segment $s_r(\tilde{G})$ with $||q - s_r(\tilde{G})||\leq \eps/(8K)\cdot ||s_r(\tilde{G}) - q'||$. 
In the point $s_r(\tilde{G})$, the function attains the value $u_0(s_r(\tilde{G}))\geq -K$, and in the endpoint $q'$, we have $u_0(q')\leq K$. 
Hence, $u_0(q')-u_0(s_r(\tilde{G}))\leq 2K$, and then 
\begin{equation*}
u_0(s_r(\tilde{G}))-u_0(q)\leq \eps/(8K) \cdot \big(u_0(q')-u_0(s_r(\tilde{G}))\big) \leq \eps/(8K) \cdot 2K = \eps/4.
\end{equation*}

As each sumand in \eqref{eq:threeterms} is at most $\eps/4$, we obtain $b(m_r(q))-u(l_{G}, q)<\eps$, or equivalently, $u(l_{G}, q)> b(m_r(q))-\eps$. 
Since $b$ is monotone increasing and concave and $m_r$ is concave, the function $b\circ m_r$ is concave. 
Hence, the conditions of Lemma~\ref{linearn} apply for $Q = S_r$, $f(p) = b\big(m_r(p)\big) - \eps$, and $f_q(p) = u(l_{\tilde{G}}, p)$ with the above $\tilde{G}$ for $q$. 
Then Lemma~\ref{linearn} provides us with the finite sequences $(q_i)$ and $(\lambda_i)$, where the $\lambda_i$ form a probability distribution (on the indices $i$). 
By using the global random variable $g$ we can simulate this distribution. 
More precisely, we can divide the unit interval into the consecutive intervals $[0,\lambda_1], [\lambda_1, \lambda_1+\lambda_2], \ldots$ of lengths $\lambda_i$, and identify each of these intervals by the unit interval. 
If $g\in [\lambda_{i-1}, \lambda_i]$ is in the $i$-th interval, then the index $i$ is chosen, and the identification of this interval with the unit interval still provides us with a uniform sample $g$ from $[0,1]$. 
Let $\tilde{G}_i$ be the graph corresponding to $q_i$. 
The mixed algorithm $l'$ then runs as follows: we pick an index $i$ according to the above law, and then run the $\PMRLA$ with fixed graph $\tilde{G}_i$, that is, the $\MRLA$ $l_{\tilde{G}_i}$. 
Note that this $\MRLA$ has the same radius $r$ as that of $l$, and for any $G\in \F$ we have 
\begin{equation*}
\E_{\omega, g}\big(\bar{v}(G, l'[\omega,g])\big) = \sum_i \lambda_i \E_{\omega,g}\big(\bar{v}(G, l_{\tilde{G}_i}[\omega,g])\big) \mathop{=}^{\eqref{udef}} \sum_i \lambda_i u\big(l_{\tilde{G}_i}, s_r(G)\big) 
\end{equation*}
\begin{equation*}
\mathop{\geq}^{\eqref{conveq}} b\Big(m_r\big(s_r(G)\big)\Big) - \eps \mathop{\ge}^{\eqref{mdef}} b\big(v^*(G)\big) - \eps. 
\end{equation*}
\end{proof}

We mention that a simplified proof to Theorem~\ref{main} can be found in \cite{Lovaszbook}. 
Here, it follows immediately from the more general Theorem~\ref{general} we have just shown. 

\begin{proof}[\emph{\textbf{Proof of Theorem \ref{main}}}]
A direct consequence of Theorem \ref{general} by putting $b(x) = x - \eps$.
\end{proof}

\section{The relations between preprocessing, mixing and randomizing}\label{sec:collapse}

We show that preprocessing and mixing are equally strong tools, and that there are no more collapses in the hierarchy other than the ones formally implied by this observation. 
The purpose of this section is to verify that the following Venn diagram is a correct illustration of the comparative expressive power of the different types of local algorithms. 
The problems separating the different sets are defined later. 
We note that all these examples separate the different types of local algorithms in the following strong sense: 
whenever a set in the Venn diagram $A$ is not contained in another set $B$, then the problem we construct is in $A$, but it is not in $B$ even with $\varepsilon$ error. 
That is, there isn't a local algorithm of type $B$ that outputs a value $\epsilon$ close to the optimum for every graph such that it might violate the allowed choices on an $\epsilon$ ratio of vertices. 
This stronger notion of not being approximable in type $B$ is relevant in graph limit theory.

Throughout this section, we slightly abuse notation in the sense that $\LA$ can refer to the set of local algorithms, and at the same time, the problems that can be solved (approximately) by some local algorithm. 
The same convention applies to all types of local algorithms. 
This convention is applied in the Venn diagram. 

On one hand, Theorem \ref{main} showed that $\MRLA$ and $\PMRLA$ are equally strong. On the other hand, it is easy to see the following theorem.

\begin{Proposition}
If a local choice problem is approximable in $\PMRLA$, then it is also approximable in $\PRLA$.
\end{Proposition}

Or more precisely,

\begin{Proposition} \label{general2}
For each correct $l_1\in \PMRLA$ there exists a correct $l_2\in \PRLA$ using the same radius $r$ and the same distribution of $\omega$ such that for each graph $G$, we have
\begin{equation*}
\E_{\omega}\Big(\bar{v}\big(G, l_2[G, \omega]\big)\Big) \ge \E_{\omega, g}\Big(\bar{v}\big(G, l_1[G, \omega, g]\big)\Big).
\end{equation*}
\end{Proposition}

\begin{proof}
If we use $l_1$ replacing $g$ with a fixed value $g_0$, then we get a $\PRLA$. 
Let us denote it by $l_1[g_0]$. 
Then 
\begin{equation} \label{average}
\E_{\omega}\Big(\bar{v}\big(G, l_1[g_0][G, \omega]\big)\Big)
\end{equation}
is a function of $g_0$. Let us choose such a $g_0$ for which \eqref{average} is at least as much as its expected value for $g_0 = g$. Let us define the $\PRLA$ $l_2 = l_1[g_0]$. Then
\begin{equation*}
\E_{\omega, g}\Big(\bar{v}\big(G, l_1[G, \omega, g]\big)\Big) = \E_g\bigg( \E_{\omega}\Big(\bar{v}\big(G, l_1[g][G, \omega]\big)\Big) \bigg)
\end{equation*}
\begin{equation*}
\le \E_{\omega}\Big(\bar{v}\big(G, l_1[g_0][G, \omega]\big)\Big) = \E_{\omega}\Big(\bar{v}\big(G, l_2[G, \omega]\big)\Big). \qedhere
\end{equation*}
\end{proof}

In particular, $\PMLA$ and $\PLA$ are also equally strong: just like before, we can cancel the randomization in Proposition~\ref{general2} by choosing a constant random variable $D$. 

This means that preprocessing and mixing are interchangeable: either one or both of them are equally strong. 
Now, we need to focus only on the strengths of mixing and randomizing. 
Obviously, every $\LA$ is a $\MLA$ and a $\RLA$, furthermore every $\MLA$ and every $\RLA$ is a $\MRLA$. 
We show that there is no further relation between the different types of local algorithms. 
Namely, we prove that $\MLA$ and $\RLA$ are incomparable: neither is stronger than the other. 
Moreover, we show that $\PMLA$ is stronger than the union of $\MLA$ and $\RLA$ by constructing a local choice problem approximable in $\PMLA$ that is not approximable in either $\MLA$ or $\RLA$. 
Finally, we prove that $\LA$ is weaker than the intersection of $\MLA$ and $\RLA$ by constructing a local choice problem approximable in both $\MLA$ and $\RLA$ that is not approximable in $\LA$. 

Some examples were known before for problems approximable in $\RLA$ but not in $\MLA$, see \cite{NgOn}.
We do not use these examples, but rather construct our own which can be used in the verification of other incomparability claims above. 
Interestingly, all examples are related to the Maximum Flow Minimum Cut Theorem. 
We only consider multisource-multitarget networks on graphs with degree at most $d$ (whose edges are directed in both ways) and all edge capacities equal to 1. 
The vertex set is partitioned into three subsets $S$, $R$, $T$, where $S$ is the set of sources and $T$ is the set of targets. 
The remaining set $R$ consist of the regular vertices, which are neither sources nor targets. 
This is easily encoded by the graph structure: sources have no loop, regular vertices have one loop, and targets have two loops. 
We call a function $f:\vec{E}(G)\rightarrow [-1,1]$ defined on the set of directed edges a \emph{fractional flow} if $f(-e)=-f(e)$ and $\sum\limits_{e \in \out(r)} f(e) = 0$ for all regular vertices $r$, where $-e$ is the reverse of the edge $e$ and $\out(r)$ denotes the set of outgoing edges from $r$. 
We call $f$ an \emph{integer flow} if $f(e)\in \{-1,0,1\}$ for all $e\in \vec{E}$. 
Normally, fractional flows are simply called flows, however we would like to emphasize the difference between the two types of flows, and introduce a certain similarity with cuts. 
The value of a flow $f$ is $\big\|f\big\| = \sum\limits_{e \in \out(S)} f(e)$, where $\out(S)$ is the set of all edges leaving $S$. 
As all edge-capacities are integers, there is a maximum fractional flow in the network that is an integer flow. 
In particular, the maximum fractional flow value coincides with the maximum integer flow value. 
A \emph{fractional cut} is a function $X: V(G)\rightarrow [0,1]$ defined on the vertex set of the graph such that $X(s)=0$ for each source $s$ and $X(t)=1$ for each target $t$. 
We call $X$ an \emph{integer cut} if $X(v)\in \{0,1\}$ for all $v\in V(G)$. 
The value of a fractional cut is $\big\|X\big\| = \sum\limits_{(a, b) \in \vec{E}(G)} c\big((a, b)\big) \max\big(0,\ {X}(b) - {X}(a)\big)$. 
As the notion of a fractional cut is more general than that of an integer cut, and the value of any fractional cut is greater than or equal to the value of any flow, the MFMC theorem implies that the minimum fractional cut value coincides with the minimum integer cut value.  
Hence, the MFMC theorem tells us that the maximum fractional flow value, the maximum integer flow value, the minimum fractional cut value and the minimum integer cut value are all equal. 

The following types of algorithms were constructed to these problems in \cite{CsE} (up to any relative error $\eps>0$): 
\begin{enumerate}
\item An $l_1\in \LA$ that finds a nearly maximum fractional flow (the $\FF$ problem). 
\item An $l_2\in \RLA$ that finds a nearly maximum integer flow (the $\IF$ problem). 
\item An $l_3\in \LA$ that finds a nearly minimum fractional cut (the $\FC$ problem). 
\item An $l_4\in \MLA$ that finds a nearly minimum integer cut (the $\IC$ problem). 
\end{enumerate}

It is easy to see that all four of these problems are indeed local choice problems. 
To be more accurate, we need to return the negative of the cut values in $\IF$ and $\IC$ to make it a maximization task: so $-\IF$ and $-\IC$ are local choice problems. 
We spell out the details for the integer cut ($\IC$); the remaining cases are left to the Reader. 

\textbf{Minimum Integer Cut Problem ($-\IC$).} Let $\F$ be the family of all graphs (or in fact, any union-closed family) such that each vertex has at most two loops and at most $d$ further edges going to other nodes. 
The problem is to find the minimum cut, that is, to partition the vertices into two sets, the source side (where $X=0$) and the target side (where $X=1$), such that the former contains all sources, the latter contains all targets, and the number of crossing edges is minimal. 
Let $C = \{0, 1, ... d, "T"\}$ be the set of choices, where $c(x) = k$ expresses that $x$ is on the source side with $k$ neighbors in the target side, and $c(x) = "T"$ means that $x$ is in the target side. 
The radius is $\delta = 1$, and $\A$ is defined as follows. 
Given the root $x$, if $c(x) \ne "T"$, then $x$ must have exactly $c(x)$ neighbors $y$ (with multiplicity) such that $c(y) = "T"$. 
If $x$ has no loop, the choice "T" is not allowed, and if $x$ has two loops, the only allowed choice is "T". 
The values are $v("T") = 0$ and $v(k) = -k$. 
This way, $\bar{v}$ expresses the negative of the relative size of the cut (whose maximum is the negative of the minimum cut, normalized by the number of vertices). 

For a graph $g$ with $n$ vertices, let $\eta(G)$ denote the vertex expansion of $G$, that is, $\eta(G)=\min\limits_{0<|S|<n/2} |N(S)\setminus S|/|S|$, where $N(S)$ denotes the set of all neighbors of ertices in $S$. 
For an $\eta>0$, let us call $G$ an $\eta$-expander if $\eta(G)\geq \eta$. 
We denote the set of all isomorphism classes of all $\eta$-expanders by $\Ex(\eta)$. 
For a graph $G$ on $2k$ vertices, let $\bisect(G)$ denote the minimum number of edges between $X \subset V(G)$ and $V(G) - X$ with $|X| = k$. 
Let $\BS(d)$ denote the limit superior of $\bisect(G) / \big|V(G)\big|$ on d-regular graphs $G$. 
The existence of expander graphs shows that $\BS(d)$ is positive. 

\begin{Proposition}\label{prop:MLARLA}
For all $l\in \RLA$ and $\eps>0$ there exists a finite graph $G$ such that the error $v^*(G) - \bar{v}(G, l[\omega])$ of $l$ for the Minimum Integer Cut Problem on the input graph $G$ is at least $\BS(d)-\eps$.
\end{Proposition}

\begin{Corollary}\label{cor:MLARLA}
The $-\IC$ problem is not approximable in $\RLA$ for $\F=\Ex(\eta)$ with small enough $\eta>0$. 
In particular, $\MLA\not\subseteq \RLA$. 
\end{Corollary}

\begin{proof}[Proof of Proposition~\ref{prop:MLARLA} and Corollary~\ref{cor:MLARLA}]
Let $r$ denote the radius the $\RLA$ uses. 
Let $G$ be a $d$-regular expander graph on $n$ vertices for a large enough $n$, so that each subset of the vertices of size $\frac{n}{2} \pm o(n)$ cuts at least $\BS(d) n - o(n)$ edges. 
Initially, we label all nodes of $G$ as regular vertices, that is, we create a network with no sources or targets. 
Let $p= \E_{\omega}\big(c^{-1}("T")/n\big)$ be the expected proportion of vertices assigned to the target side by the algorithm $l$ executed on this network. 
If $p>1/2$, let us change the status of nodes in $G$ one by one to source; if $p<1/2$, then we change the status of nodes to target, instead. 
Changing the status for one node only affects the choices in the $r$-neighborhood of the node, which contains at most $(d + 1)^r$ vertices. 
If we change the status of all nodes, then the expected ratio $\E_{\omega}\big(c^{-1}("T")/n\big)$ decreases to 0 or increases to 1, depending on whether $p>1/2$ or $p<1/2$. 
Therefore, we can stop the procedure at a point when this expected proportion is $\frac{1}{2} \pm o(1)$. 
The network $G$ obtained this way has no targets or no sources, thus its minimum cut value is $v^*(G) = 0$. 

The choice at each node is independent from all but at most $(d + 1)^{2r}$ nodes, thus 
\begin{equation*}
\Var\Big(\frac{c^{-1}("T")}{n}\Big) = \frac{1}{n^2} \sum_{x, y \in V(G)} \Cov\Big(c(x) = "T", c(y) = "T"\Big) \le \frac{1}{n^2} \sum_{x \in V(G)} (d + 1)^{2r} = \frac{(d + 1)^{2r}}{n}. 
\end{equation*}
Hence, $c^{-1}("T") / n = \frac{1}{2} + o(1)$ with high probability, making the expected size of the cut $\BS(d) n - o(n)$.

Applying this observation to the Minimum Cut Problem as defined above, $v^*(G) = 0$, while $\bar{v}(G, c)$ can be arbitrarily close to $-\BS(d)$, yielding an error $v^*(G) - \bar{v}(G, c)$ arbitrarily close to $\BS(d)$.
\end{proof} 

Hence, the $-\IC$ problem is approximable in $\MLA$ but not in $\RLA$ for $\F=\Ex(\eta)$ with small enough $\eta>0$.
We show that $\IF$ is not approximable in $\MLA$ (but as we mentioned earlier, it is approximable in $\RLA$) for $\F=\Ex(\eta)$ with small enough $\eta>0$. 

\begin{Proposition}
Let $d,k\in \mathbb{N}$, $d\geq 2$, and let $A$ and $B$ be two disjoint sets with size $k$. 
Let us produce a network $G$ from a $d$-regular bipartite graph $G_0$ with bipartition $A,B$ by linking a source with an (undirected) edge to each vertex of $A$ and a target to each vertex in $B$. 
Let $l\in \MLA$ have radius $r$ be a correct algorithm for the $\IF$ problem for $\F=\Ex(\eta)$ with small enough $\eta>0$. 
If the girth of $G_0$ is at least $2r+2$, then $l$ outputs the all-zero flow on $G_0$, and thus $v^*(G_0) - \bar{v}(G_0, l[g])\geq 1/4$. 
\end{Proposition}

\begin{proof}
First of all, note that such graphs exist in $\Ex(\eta)$: a random $d$-regular bipartite graph has large girth, and it is an $\eta$-expander for some fixed $\eta$ (that depends on $d$). 
The minimum flow clearly has value $k$. 
In fact, such a fractional flow can easily be constructed by assigning 1 to edges starting at a source and to edges ending at a target, and  $1/d$ to edges between $A$ and $B$. 
An $\MLA$ might find this fractional flow. 
However, we show that $\MLA$s struggle to find integer flows. 

We make use of the symmetry that cannot be broken by an $\MLA$, since we are not allowed to assign random seeds to the vertices. 
Namely, as the girth of $G_0$ is at least $2r+2$, the $r$-neighborhood of any edge in $G_0$ is the same: the regular vertices form a tree with a root edge of depth $r$ such that every vertex that is not a leaf has degree $d$, and each regular vertex has an additional neighbor that is a source or a target, in an alternating fashion. 
Hence, an $\MLA$ must assign the same value to any edge in $G_0$. 
This common value cannot be 1, as only 1 unit of flow can enter a vertex in $A$, thus it violates the conditions of a flow if $d$ units leave such a vertex. 
Therefore, the common value assigned to every edge in $G_0$ is 0. 
Thus the optimal (integer) flow ratio is $k/4k=1/4$, and any correct $\MLA$ outputs the all-zero flow, with value 0. 
\end{proof}

\begin{Corollary}
$\RLA\not\subseteq \MLA$
\end{Corollary}

Now we have shown that neither one of $\MLA$ and $\RLA$ is stronger than the other. 
Finally, we show that $\MRLA$ is strictly stronger than the ``union'' of $\MLA$ and $\RLA$, and that $\LA$ is strictly weaker than the ``intersection'' of $\MLA$ and $\RLA$.

\begin{Proposition}
$\MLA\cup \RLA \subsetneq \MRLA$
\end{Proposition}
\begin{proof}
The idea is to ``add up'' the problems of $\IF$ and $-\IC$ to obtain the problem $\IF-\IC$. 
(Here, we can choose $\F$ as the set of isomorphism classes of all finite d-regular graphs.) 
That is, the local algorithm should construct both an integer flow and a fractional cut, and the value is the difference of the flow and cut values. 
Hence, a near-optimal solution must consist of a near-optimal integer flow and a near-optimal integer cut, as the two errors add up. 
Then this local choice problem is approximable neither in $\MLA$ nor in $\RLA$, but it is approximable in $\MRLA$. 
\end{proof}

The other claim requires a more elaborate argument. 

\begin{Proposition}\label{prop:meet}
$\LA \subsetneq \MLA\cap \RLA$
\end{Proposition}
\begin{proof}
The idea is to dualize the previous construction. 
That is, rather than solving two local choice problems simultaneously, we would like to define a local choice problem that can choose between solving two local choice problems: one that is approximable in $\RLA$ but not in $\MLA$, and one that is approximable in $\MLA$ but not in $\RLA$. 
Then we can find an approximate solution if we can apply either an $\RLA$ or an $\MLA$, but an $\LA$ is insufficient for the task. 
There are two major difficulties with this strategy (in particular we cannot use the problems $\IF$ and $-\IC$ again). 
We need to make sure that both problems have the same optimum, otherwise the choice is not arbitrary: the algorithm would be forced to approximate the one with the larger optimum. 
Moreover, we need to make sure that the algorithm indeed chooses one out of the two problems, rather than outputting a mixed solution. 
This is achieved by penalizing such mixtures in local structures. 

This time, $\F=\Ex(\eta)$ for a small enough $\eta>0$. 
We start by defining the difference of two pairs of flow-related problems: $\IF-\FC$ and $\FF-\IC$. 
For example, $\IF-\FC$ is the local choice problem where we need to find an integer flow and a fractional cut simultaneously, and the value is the difference of the flow and cut values. 
Clearly, the maximum value is 0, and a near-optimal solution can be obtained if and only if we can find a near-optimal integer flow and a near-optimal fractional cut. 
Hence, this problem is approximable in $\RLA$ but not in $\MLA$. 
The other problem $\FF-\IC$ is defined analogously. 
Once again, the maximum value is 0, and a near-optimal solution can be obtained if and only if we can find a near-optimal fractional flow and a near-optimal integer cut. 
Hence, this problem is approximable in $\MLA$ but not in $\RLA$. 

Now that we have the local choice problems $\IF-\FC$ and $\FF-\IC$, we define the disjunction $(\IF-\FC)\vee (\FF-\IC)$. 
Informally, we want to make sure that a near-optimal solver picks $\IF-\FC$ or $\FF-\IC$ and solves it. 
So we take the union of all colorings (allowed pairs) for $\IF-\FC$, and $\FF-\IC$, adding an extra red color for local structures in the former and an extra blue color to those of the latter. 
We are allowed to ignore one of the problems. 
That is, a correct coloring consists of red, blue and mixed red and blue (purple) local structures. 
Given a vertex $u$ that is colored red, if $u$ or any neighbor of $u$ is colored blue, then on top of the usual values, we add a penalty of $-4/\min(\eta,1)$. 
The same happens if a blue vertex has a red neighbor. 
The final score is the sum of all flow and negative cut values together with the penalties.  
Clearly, the maximum attainable value is 0: this is the optimum of $\IF-\FC$ and $\FF-\IC$, and the penalties only decrease the value. 
Furthermore, 0 is indeed attainable in two trivial ways: we solve the $\IF-\FC$ problem (by providing a maximum integer flow and a minimum fractional cut) and ignore the other problem, or we solve the $\FF-\IC$ problem and ignore the other problem. 
In particular, the problem $(\IF-\FC)\vee (\FF-\IC)$ is approximable in $\RLA$ and in $\MLA$, as well. 

We show that the problem $(\IF-\FC)\vee (\FF-\IC)$ is not approximable in $\LA$. 
If in a correct solution, the fraction of purple vertices is more than $\eps/2$, then their overall penalty amounts to a negative number below $-2\eps$. 
Hence, such a solution cannot be nearly optimal, as the value of a nearly optimal solution is at least $-\eps$, and the ``flow minus cut'' type contribution in the value is non-positive. 
Let $c$ be the color that occurs less frequently out of red and blue. 
Then less than half of the vertices is colored only by $c$ and not the other color. 
If the fraction of such vertices is more than $\eps/2$, then they have a neighborhood containing at least an $\eps\eta/2$ fraction of the vertices due to the expansion in the graph, all of whose elements are colored by the other color (different from $c$). 
The overall penalty given to this neighborhood amounts to a negative number below $-2\eps$. 
Hence, such a solution cannot be nearly optimal. 
Thus in a nearly optimal solution, there is a dominant color, and the other color occurs in an at most $\eps$ fraction of the vertices. 
There is a $\varrho>0$ such that if $k$ is large enough and $G_0$ is a random bipartite graph, then with high probability, the girth of $G_0$ is at least $2r+2$, and every $r$-neighborhood that can occur in such a construction occurs with frequency at least $\varrho$. (Set $\varrho$ as half of the minimum expected frequency of an $r$-neighborhood in a large random construction.) 
If $\eps<\varrho$, then the above observation shows that a nearly optimal solution is monochromatic: it is either a correct solution of $\IF-\FC$ with the $\FF-\IC$ part ignored, or the other way around. 
Whichever the case, it cannot be achieved by an $\LA$, since neither $\IF-\FC$ nor $\FF-\IC$ is approximable in $\LA$. 
\end{proof}

We note that we can define a similar counterexample where $\F$ consists of all finite graphs with degree at most $d$. 
We simply choose the penalty as $-4/\min(\eta,1)$ for some $\eta$ such that a sequence of $\eta$-expanders exists. 
Moreover, we could also pick $\F$ to be the class of all finite graphs with degree at most $d$ throughout the propositions in this section, wherever we chose $\F=\Ex(\eta)$.

Finally, we introduce an intriguing variant of preprocessing. 
Let us call a local algorithm \emph{weakly statistical} if for any $\eps>0$ there is a fixed radius $r$ such that $c(x)$ depends on $B_r(x)$ and a random vector that is closer in $L_2$ distance to $s_r(G)$ than $\eps$ with probability at least $1-\eps$. 
That is, we are allowed to draw a large sample of $r$-neighborhoods from the graph to find an empirical neighborhood statistics. 
One would expect that this is all the information about the isomorphism type of the graph that might be relevant for a local choice problem. 

\begin{Question}\label{que:W}
Is $\WRLA=\MRLA$, where $W$ stands for weakly statistical?
\end{Question}

As before, a positive answer would imply $\WLA=\MLA$, as well. 
We note that Lemma~\ref{linear} yields a positive answer to Question~\ref{que:W} if the valuation function $v$ of the local choice problem is bounded. 
So the problem is only interesting for local choice problems where the valuation function is not bounded from below. 

Some applications of local algorithms (e.g. factor of iid processes) are less strict about allowed colorings.
As already indicated at the beginning of this section, we say that an algorithm $l$ is $\varepsilon$-correct if for all graphs $G$ we have $\mathbb{P}_{\omega, g, x}\left(\big(B_{\delta}(G, x), l[G, \omega, g]\upharpoonright_{V(B_{\delta}(G, x))}\big) \in \A\right) > 1 - \eps.$
We call a local choice problem weakly approximable in a given type of local algorithms if for all $\varepsilon>0$ there is an algorithm $l$ of that type that produces an $\varepsilon$-correct solution $\varepsilon$-close to the optimum.
Each of our examples in this section was such that whenever we claimed it was not approximable in a given type, it was not weakly approximable either. 
By dropping this additional requirement one could somewhat simplify the examples, especially the one in the proof of Proposition~\ref{prop:meet}. 

\section{Acknowledgements}

The author is thankful to L\'aszl\'o Lov\'asz, G\'abor Elek, G\'abor Lippner and Mikl\'os Ab\'ert for their guidance, and to Andr\'as Pongr\'acz for the useful comments on the presentation of the paper. 
The author was supported by the NRDI grant KKP~138270.

\bibliographystyle{plain}
\bibliography{refs}
\end{document}